\documentclass[a4paper]{article}
\usepackage{amsmath,amsthm,amssymb,enumerate,hyperref,bm,color,eurosym}
\usepackage[legacycolonsymbols]{mathtools}
\providecommand\given{}
\newcommand\SetSymbol[1][]{%
  \nonscript\:#1\vert
  \allowbreak
  \nonscript\:
  \mathopen{}}
\DeclarePairedDelimiterX\Set[1]{\{}{\}}{%
  \renewcommand\given{\SetSymbol[\delimsize]}
  #1}

\DeclarePairedDelimiterXPP\pospart[1]{}{(}{)}{^+}{#1}
\DeclarePairedDelimiterXPP\negpart[1]{}{(}{)}{^-}{#1}

\newcommand\R{\mathbb{R}}

\newcommand\N{\mathbb{N}}
\newcommand\PP{\mathcal{P}}

\newcommand\GG{\mathcal{G}}

\newtheorem{thm}{Theorem}[section]
\newtheorem{prop}[thm]{Proposition}

\theoremstyle{definition}
\newtheorem{defn}[thm]{Definition}
\newtheorem{ex}[thm]{Example}
\newtheorem{rem}[thm]{Remark}

\DeclareMathOperator*{\vertices}{vert}

\DeclareMathOperator*{\gr}{gr}
\DeclareMathOperator*{\cl}{cl}
\DeclareMathOperator*{\conv}{conv}

\DeclareMathOperator*{\Min}{Min}
\DeclareMathOperator*{\EVPI}{EVPI}

\makeatletter
\newcommand{\leqnomode}{\tagsleft@true\let\veqno\@@leqno}
\newcommand{\reqnomode}{\tagsleft@false\let\veqno\@@eqno}
\makeatother

\usepackage{booktabs}

\newcommand{\lel}{\preccurlyeq}

\definecolor{color1}{rgb}{0,.2,.8}
\definecolor{color2}{rgb}{1,.2,0}
\definecolor{color3}{rgb}{.2,.7,.6}

\title{Multi-objective stochastic linear programming with recourse and flexible decision making}
\author{Andreas H. Hamel\thanks{Free University of Bozen-Bolzano, Italy, andreas.hamel@unibz.it} \and Andreas L{\"o}hne\thanks{Friedrich Schiller University Jena, Germany, andreas.loehne@uni-jena.de}}

\begin{document}
\maketitle

\begin{abstract} \noindent
Optimal inventory leads to stochastic optimization problems where deterministic delivery decisions have to be made in advance of stochastic demand realizations. Similarly, risk deposits have to be given before the random outcomes of investments are known. In this paper, multi-criteria versions of such stochastic recourse problems are studied.

In addition to traditional concepts like Pareto-optimality, a decision maker for the multi-criteria problem may have a preference for greater flexibility in the second stage decision. This idea
leads to a first stage optimization problem with a set-valued objective instead of a mere multi-criteria one. Under linearity assumptions, this problem becomes a polyhedral convex set optimization problem instead of a multi-objective linear program. 

Solution concepts for multi-objective/set-valued recourse problems are given as well as deterministic surrogates of the stochastic problem such as its deterministic equivalent, the so-called wait-and-see problem and the expected-value problem for the multi-objective case. Managerial decision making guidelines are obtained based on set optimization methods and the preference-for-flexibility approach: choose the deterministic first stage variable such that a maximum of flexibility is combined with a guarantee for Pareto minimal objective values.
Two major examples illustrate the findings, a multi-objective newsvendor problem with an additional health/sustainability objective and a risk compensation problem where the availability of more than one asset for risk compensation, e.g., several currencies, leads to multiple objectives.



\end{abstract}

\section{Introduction}
\label{SecIntroduction}

Stochastic optimization problems with recourse \cite{WalkupWets67SIAppMath} are characterized by a two stage decision making procedure: later stage decisions depend on uncertain outcomes of the earlier stage. Typically, the first stage decision variable is deterministic and has to be chosen before the random data of the problem are realized. The objective value for the second stage decision depends on both the realization of the random data and the choice of the first-stage variables. The goal is to optimize this objective value's expectation subject to all feasible first-stage variables.

Such problems with a single (real-valued) objective were widely studied and applied to many management problems prominent among them the optimal inventory policy problem. This paper deals with the multi-objective stochastic linear programming problem with recourse. Our primary objective is to delineate a decision-making procedure for such problems that is practical and tractable.

We argue that classical decision-making techniques from multi-objective optimization (e.g., just look for (elements of) the Pareto frontier) exhibit certain limitations when naively applied to recourse problems. Indeed, a decision maker might be interested in keeping second-stage decisions as flexible as possible in order to be able to react to realizations of the random data. This potential goal adds another layer to the original objectives and deserves special treatment. It is also a special feature of multi-objective problems since it does not produce anything new in the single-objective case.

To handle this preference for flexibility, we introduce a polyhedral convex set optimization problem as an alternative to a multi-objective linear program reformulation of the original stochastic problem. We illustrate that solving this set optimization problem provides the decision maker a set of alternatives, being Pareto optimal with respect to the original objectives of the problem, and simultaneously maintaining a maximum of flexibility for decisions in the second stage. A multi-objective generalization of the newsvendor problem serves as a showcase application.

Moreover, we establish surrogate problems parallel to the single-objective case. In particular, the multi-objective wait-and-see problem is formulated which admits to address the question of the value of information: how much could the optimal value of the recourse problem be improved if the first-stage decision could be made after the random data has been realized? In the multi-objective case, the optimal values turn out to be sets which can be interpreted as infima in an appropriate complete lattice of sets. 

Similarly, the multi-objective expected value problem is given which is obtained by replacing the random data of the model by their expectations. Following the same idea as before, a corresponding polyhedral convex set optimization problem is given which can be used to define the value of the stochastic solution by comparing sets.

The methodology for this paper borrows from two sources. The first one is the dynamic structure of decision making with respect to preferences for flexibility which are studied in economics, but rarely involved in managerial decision making: at a first stage/earlier time, a decision among sets of alternatives is made which keeps flexibility for choosing a particular alternative from the chosen set only later \cite{Kreps79Economet}. The second one is the complete lattice approach to set optimization \cite{HamelEtAl15Incoll} which provides a comprehensive framework for comparing sets and solving set-valued optimization problems as well as an extension of much of the theory from single- and multi-criteria to set-valued optimization problems.

The rest of the paper is organized as follows. Section \ref{SecLiterature} is a literature review. The notation for the two-stage, multi-objective stochastic linear program with recourse is introduced in Section \ref{SecLRP} as well as some (solution) concepts from multi-objective optimization in such a way that the link to optimization with respect to preferences for flexibility can be made later on; moreover, a bi-objective version of the classic newsvendor problem is discussed as a motivating example. This example, which has an additional sustainability objective, is also used in the following sections to illustrate the novel concepts and results in the paper.

In Section \ref{SecSOandPfF}, a set optimization framework with appropriate solution concepts for the decision making problem with respect to preferences for flexibility is presented. It is specialized to the multi-stage multi-objective problem in Section \ref{SecMDMwithPF}.
Section \ref{SecWSP} is devoted to the wait-and-see problem which admits to define the expected value of perfect information for the multi-objective problem, while Section \ref{SecEVP} provides the expected-value problem. In Section \ref{section:ex_risk}, an application to risk management is discussed involving a market model with transaction costs and multiple eligible portfolios. In the Appendix, proofs of a few mathematical facts are collected.

\section{Literature review}
\label{SecLiterature}

Multiple stage stochastic linear programs were independently introduced by George Dantzig \cite{Dantzig55MS} and Martin Beale \cite{Beale55} in 1955. This marked the beginning of a major research area. An important instance of such problems is the optimal inventory policy problem \cite{ArrowHarrisMarschak51Economet}, see also \cite{GeigerLodree03Proc} and the references therein. An alternative approach can be found in the recent \cite{BertsimasEtAl23MS} following a robust optimization idea. An extension to nonlinear problems was already given in \cite{MangasarianRosen64OR, WalkupWets67SIAppMath}.

If the involved functions are linear and the random coefficients are finitely distributed, the recourse problem can be expressed equivalently as an ordinary linear program as already done in \cite{Beale55, Dantzig55MS}. This usually leads to large problem instances since this reformulation requires a duplication of the second-stage variables and the corresponding constraints according to the number of scenarios of the distribution. 

Alternatively, surrogate problems have been formulated among them the wait-and-see problem \cite{AvrielWilliams70OR}. The difference between the optimal objective values of the wait-and-see problem and the recourse problem is called expected value of perfect information, see \cite{AvrielWilliams70OR} and the references therein. 

The expected-value problem is a simplification of the recourse problem avoiding the duplication of variables and constraints. The expected worsening of this outcome is referred to as value of the stochastic solution \cite{Birge82MP}.

A comprehensive examination of the relationships between the recourse problem, the wait-and-see problem and the expected-value problem for single-criteria problems can be found in \cite{Birge82MP}.   

The literature on multi-objective recourse problems is sparse (compare the survey \cite{GutjahrPichler16AOR}) but potential applications seem to be very interesting. As examples, we point at one to recommender systems \cite{GhanemLeitnerJannach22ECRA} and at another to disaster relief strategies \cite{GendreauGutjahrRath16ITOR}. The recent \cite{GuptaHunter25} discusses linear two-stage problems on general probability spaces and points towards applications in many areas. Therein, it is stated: `Ultimately, the optimization problem $\ldots$ is a {\em stochastic set optimization} problem' \cite[Section 1.1]{GuptaHunter25}, but the solution concept is different from the one proposed in this paper.

The idea of preferences for flexibility based on order relations for subsets of an ordered set, nowadays called set relations, was first discussed in \cite{Kreps79Economet} (see \cite{HamelLoehne20MMOR} for a more recent presentation). The two-stage character of stochastic recourse problems fits perfectly into the framework of set relations since also the latter has a dynamic feature: make a decision for a set of alternatives now and choose a particular alternative from the chosen set only later.

However, Kreps \cite{Kreps79Economet} assumed a complete base preference, i.e., every pair of alternatives is comparable. This assumption was dropped in later publications \cite{ArlegiNieto01MSS, DeMouraRiella21TD} and is also not made in the present paper due to the multi-objective nature of the problem. Applications in economics include the comparison and ranking of so-called opportunity sets \cite{ArlegiNieto01SCW, BarberaGrodal11JME}.

To the best of our knowledge, the preference-for-flexibility approach was never really applied to managerial problems which probably is due to the lack of an optimization theory for set-valued functions. This gap was filled by the complete lattice approach to vector and set optimization which is surveyd in \cite{HamelEtAl15Incoll} providing tools for the development of the present paper.

\section{The multi-objective linear recourse problem}
\label{SecLRP}

First, some notation is fixed. Let $(\Omega, p)$ be a {\em finite probability space}, that is, a nonempty finite set $\Omega\coloneqq\Set{\omega_1, \dots,\omega_N}$ and a function $p \colon \Omega \to (0,1]$ with $\sum_{\omega \in \Omega} p(\omega) = 1$. 

For an arbitrary set $R$ we denote by $R^\Omega$ the family of all functions from $\Omega$ to $R$. This leads, for instance, to the notation $p \in (0,1]^\Omega$. Since $\Omega$ is a finite set of cardinality $N$, $R^\Omega$ can be identified with $R^N \coloneqq R \times \dots \times R$. For $r \in R^\Omega$, we set $r_\omega \coloneqq r(\omega)$. For a {\em random matrix} $\bm{M} \colon \Omega \to \R^{m\times n}$ (for $n=1$ also called {\em random vector}) or a {\em random set}  $\bm{S} \colon \Omega \to 2^{(\R^d)}$ we write
$$ \bm{M}= \begin{pmatrix}
	M(\omega_1)\\ \vdots \\ M(\omega_N)
\end{pmatrix} \in (\R^{m\times n})^\Omega \quad \text{ and } \bm{S}= \begin{pmatrix}
	S(\omega_1)\\ \vdots \\ S(\omega_N)
\end{pmatrix} \in (2^{(\R^d)})^\Omega,$$
so $\bm{M}$ can be seen as a ``vector of matrices'' and $\bm{S}$ is a ``vector of sets''.

Let $\bm{M} \in (\R^{m \times n})^\Omega$ and $\bm{X} \in (\R^{n \times k})^\Omega$ be random matrices. Then we define the product
$$ \bm{M}\bm{X} \coloneqq \begin{pmatrix}
	M(\omega_1)X(\omega_1)\\
	\vdots \\
	M(\omega_N)X(\omega_N)\\
\end{pmatrix} \in (\R^{m \times k})^\Omega,$$
which refers to a scenario-wise matrix multiplication.
We also use the notation $v^\omega\coloneqq v(\omega)$ for vectors and $M_\omega\coloneqq M(\omega)$ for all other objects (because the lower index for vectors is reserved for its components).

The {\em expectation} of a random set $\bm{S}$ is the set $E[\bm{S}]\coloneqq \sum_{\omega\in \Omega} p_\omega S_\omega$, where the sum of sets in $\R^d$ is defined by the usual Minkowski addition. 
If $\bm{S}$ is a random set such that $S_\omega$ is a polyhedron (always assumed to be convex here), $\bm{S}$ is also called a {\em random polyhedron}. Since every polyhedron has an inequality representation, a random polyhedron $\bm{P}$ can be expressed by a random matrix $\bm{B} \in (\R^{m\times n})^\Omega$ and a random vector $\bm{c} \in (\R^m)^\Omega$ as
$$ \bm{P} = \begin{pmatrix}
	\Set{x \in \R^n \given B(\omega_1) x \geq c(\omega_1)}\\
	\vdots \\
	\Set{x \in \R^n \given B(\omega_N) x \geq c(\omega_N)}
\end{pmatrix} \in (2^{(\R^n)})^\Omega.$$
For short we write
$$  \bm{P} \coloneqq \Set{\bm{x} \in (\R^n)^\Omega \given \bm{B} \bm{x} \geq \bm{c}}. $$
A set-valued mapping $F:\R^n \rightrightarrows \R^d$ is called \emph{polyhedral convex} if its graph 
$$ \gr F \coloneqq \Set{(x,y)\in \R^n \times \R^d \given y \in F(x)}$$
is a convex polyhedron.
A {\em random polyhedral convex set-valued function} $\bm{F}: \R^n \to (2^{\R^d})^\Omega$ is an $N$-tuple of polyhedral convex set-valued mappings $F(\cdot,\omega):\R^n \rightrightarrows \R^d$ defined as
$$ \bm{F}(x) \coloneqq  \begin{pmatrix}
	F(x,\omega_1)\\
	\vdots \\
	F(x,\omega_N)
\end{pmatrix} \in (2^{\R^d})^\Omega.$$
\begin{rem}
	A ``random set-valued mapping'' $\bm{F} \colon \R^n \rightrightarrows (\R^d)^\Omega$ is different from a random set-valued function $\bm{F} \colon \R^n \to (2^{\R^d})^\Omega$. The reason is that for $N = |\Omega| \geq 2$,
	$$ 2^{(\R^d)^\Omega} \neq (2^{\R^d})^\Omega.$$
	For instance, for $d=1$, $N=2$ and the set $S$ the line segment between the unit vectors in $\R^2$, one has $S \in 2^{\R \times \R}$ but $S \not\in  2^\R \times 2^\R$.
\end{rem}

The expectation of a random polyhedral convex set-valued function is the polyhedral convex set-valued mapping
$$ E[\bm{F}]:\R^n \rightrightarrows \R^d,\quad E[\bm{F}](x)\coloneqq \sum_{\omega\in \Omega} p(\omega)F(x,\omega).$$

A  {\em multi-objective stochastic linear program with recourse}, abbreviated as a {\em recourse problem} here, is the problem to
\leqnomode
\begin{gather*}\tag{RP}\label{spr}
\text{minimize}\;\;  Cx +  E[\bm{Q}\bm{y}] \quad\text{s.t.} \quad
\left\{ 
  \begin{array}{rcl}
  	Ax&=&b\\
  	\bm{T} x + \bm{W}\bm{y} &=& \bm{u}\\
  	x, \bm{y} &\geq& 0\\
	 x \in \R^n,&&\hspace{-.8cm}\bm{y} \in (\R^m)^\Omega\text{.}
  \end{array}\right. 
\end{gather*}
\reqnomode
where the decision variable $x \in \R^n$ and the data $C \in \R^{d\times n}$, $A \in \R^{k \times n}$, $b \in \R^k$ are deterministic, while the decision variable $\bm{y} \in (\R^m)^\Omega$ and the data $\bm{Q} \in (\R^{d \times m})^\Omega$, $\bm{T} \in (\R^{\ell \times n})^\Omega$, $\bm{W} \in (\R^{\ell \times m})^\Omega$, $\bm{u} \in (\R^\ell)^\Omega$ are random. 
Assumption~(A1) below clarifies the order relation for the objective values and the meaning of ``minimize'', while one of the goals of this paper is to provide a precise interpretation of the minimization procedure.

This recourse problem models a two-stage decision process.  The first-stage variables $x_1, \dots, x_n$ must be fixed before the random data is observed. Therefore, the expectation of the objective value is optimized. The second-stage decision is made after the random coefficients are realized. It not only depends on the first-stage decision $x$ but also on the realization of random data. Therefore, the vector $\bm{y}$ of second-stage variables is random.

Problem \eqref{spr} can be seen as a multi-objective linear program which can be expressed equivalently as
\leqnomode
\begin{gather*}\tag{RP'}\label{spr1}
  \text{minimize}\;\;  C x + \sum_{i=1}^N p_i Q_i y^i \hspace{6cm}\\ 
  \quad\text{s.t.} 
  \left\{ 
  \arraycolsep=1.0pt
  \begin{array}{rlrlcrlrl}
  	A     &x&       &   &        &       &   & =& b \\
  	T_1   &x& + W_1 &y^1&        &       &   & =& u^1\\
	\vdots& &       &   & \ddots &       &   &  &    \\
	T_N   &x&       &   &        & + W_N &y^N& =& u^N\\
       	\multicolumn{7}{r}{x,\,y^1,\dots,y^N}&\geq&0 \\
	&&&&&& x \in \R^n, \, y^1,\dots,y^N & \in & \R^m \text{,} \nonumber
  \end{array}
  \right.
\end{gather*}
\reqnomode
where we set $Q_i \coloneqq Q(\omega_i) \in \R^{d \times m}$, $T_i \coloneqq T(\omega_i) \in \R^{\ell \times n}$, $W_i \coloneqq W(\omega_i) \in \R^{\ell \times n}$, $u^i \coloneqq u(\omega_i) \in \R^\ell$ and $p_i \coloneqq p(\omega_i) \in [0,1]$.
The primary goal is to optimally determine the vector $x$, representing the optimal decision in the first stage. 
For each first-stage decision $x$ and each scenario $\omega \in \Omega$, the optimal second-stage decision $y= y(x,\omega)$ is obtained from the multiple objective linear program 
\leqnomode
\begin{gather*}\tag{RP$_2(x,\omega)$}\label{spr2xo}
\text{minimize}\;\;  Cx +  Q(\omega) y \quad\text{s.t.} \quad
\left\{ 
  \begin{array}{rcl}
  	 W(\omega) y  &=& u(\omega) - T(\omega) x\\
  	y &\geq& 0	\\
	&&\hspace{-.8cm}y \in \R^m \text{.}
  \end{array}\right. 
\end{gather*}
\reqnomode

In order to be able to explain several concepts from multi-objective linear programming in a transparent way as described, for example, in \cite{Ehrgott05Book, Loe11}, we introduce a more compact version of Problem \eqref{spr1}. Let $S$ denote the feasible set, i.e., the set of all $z=(x,y^1,\dots,y^N) \in \R^{n+Nm}$ satisfying the constraints of \eqref{spr1} and $P \coloneqq (C,p_1 Q_1,\dots,p_N Q_N) \in \R^{d \times (n+Nm)}$ the objective matrix of \eqref{spr1}. Problem \eqref{spr1} can now be re-written as
\leqnomode
\begin{equation*}\tag{RP''}\label{spr2}
  { \text{minimize}} \quad P z \quad \text{s.t.} \quad z \in S.
\end{equation*}
\reqnomode
In this paper, we make two further assumptions for the sake of technical transparency. The first one is
\begin{description}
\item{(A1)}
 the outcomes of Problem \eqref{spr1}, i.e., the objective values, are ordered component-wise: the ordering cone is $\R^d_+$.
\end{description}

The {\em upper image} of \eqref{spr2} is the convex polyhedron
\begin{equation}\label{eq:upper_image_spr2}
 \PP \coloneqq \bigcup_{z \in S} \Set{Pz} + \R^d_+.
\end{equation}
It conveys the set of all possible outcomes $Pz$ for $z \in S$ together with all points being dominated by such outcomes with respect to the component-wise partial ordering $\leq_{\R^d_+}$. Thus, the upper image includes the possible outcomes together with all hypothetical outcomes being worse than these. The upper image is related to the free disposal condition which is basic in economics, see \cite[p. 131]{MasColellEtAl95Book}, for example. The set  $\PP$ satisfies $\PP + \R^d_+ = \PP$, i.e., it satisfies a ``positive" free disposal condition due to the fact that minimization is the goal. With the decision making problem \eqref{spr1} in view, knowing $\PP$ means knowing all one's options.

The second assumption is
\begin{description}\label{pseudo_label}
\item{(A2)} the recourse problem \eqref{spr1}, \eqref{spr2} is \emph{bounded}, which means
\begin{equation}\label{eq:bounded}
  \PP = \conv (\vertices \PP) + \R^d_+\text{,}
\end{equation}
where $\vertices \PP$ denotes the set of vertices of $\PP$ and $\conv S$ the convex hull of a set $S \subseteq \R^d$.
\end{description}

A point $\bar v \in \PP$ is called {\em minimal} in $\PP$ if there is no $v \in \PP$ such that $v \in \Set{\bar v} - \R^d_+\setminus\Set{0}$: if $v$ is better ($=$smaller) than $\bar v$ in one component then it must be worse (larger) in at least one other. 

Clearly, knowing all one's options, i.e., the set $\PP$, is not enough---the decision maker wants to find ``best solutions.'' Candidates for such solutions are points $\bar z \in S$ such that $P \bar z$ is minimal in $\PP$. It is known that the vertices of the upper image $\PP$ are minimal \cite[Corollary 4.67]{Loe11}. So, a simple approach for the decision maker could be to choose a vertex of $\PP$ but it is not a priori clear which one to pick. Moreover, in general, there exist many more minimal points of $\PP$. 

On the other hand, if (A2) is satisfied, every outcome is dominated by a convex combination of vertices, which follows from \eqref{eq:bounded}. In particular, every minimal point of $\PP$ is a convex combination of vertices of $\PP$ (while a convex combination of vertices is not necessarily a minimal point). A feasible point $\bar z \in S$ such that $P\bar z$ is minimal in $\PP$ is called a {\em minimizer} of \eqref{spr1} (or {\em Pareto-minimal} or {\em efficient}, see \cite[Definition 2.1]{Ehrgott05Book}). The facts that all vertices of $\PP$ are minimal and, together, they generate the upper image via \eqref{eq:bounded} motivate the following definition.

\begin{defn}
\label{DefSolution} 
A finite set $\Set{z^1,\dots,z^k} \subseteq S$ of minimizers is called {\em solution} of the multi-objective linear program \eqref{spr2} if $\vertices \PP \subseteq \Set{Pz^1,\dots,Pz^k}$.
\end{defn}

This solution concept subsumes the idea of ``choosing best alternatives'' as well as ``knowing all one's options:" each minimal element of $\PP$ can be generated via convex combinations of $\Set{Pz^1,\dots, Pz^k}$ and $\PP$ can be generated via
\[
\PP = \conv \Set{Pz^1,\dots, Pz^k} + \R^d_+\text{.}
\]
Below, we will show that this solution concepts fits with preferences for flexibility as well as set optimization approaches \cite{HamelEtAl15Incoll, Loe11}. 
We redirect our focus to the recourse problem.

\begin{prop}\label{prop:31}
Let $\bar z \coloneqq (\bar x,\,\bar y^1,\dots,\bar y^N) \in S$ be a minimizer of \eqref{spr1}. Then for all $i \in \Set{1,\dots,N}$, $\bar y^i$ is a minimizer of \rm{(RP$_2(\bar x,\omega_i)$)}.
\end{prop}
In view of Proposition \ref{prop:31}, it appears to be sufficient to solve the multiple objective linear program \eqref{spr1} and let the decision maker act in the first stage only. However, the decision maker's preferences (in the sense of choosing suitable weights for the objectives) can change after the first-stage decision has been made. This means the decision maker might have a {\em preference for flexibility} \cite{Kreps79Economet}, i.e., the desire to have as many as possible alternative options for her second-stage decisions. Let us illustrate this idea by an example.

\begin{ex}\label{ex:1} (Multi-objective newsvendor problem I)\\
We consider a newsvendor who sells $n$ types of newspapers. In the morning she purchases $x=(x_1,\dots,x_n)^T \in \R^n$ of each type of newspaper in order to sell them during the day. The demand of the various types of newspapers is random and has been observed over $N$ days $\Omega=\Set{\omega_1,\dots,\omega_N}$ in the past. Let $d_i(\omega)$ be the demand of newspaper $i \in \Set{1,\dots,n}$ at day $\omega \in \Omega$ during a maximal daily working time (say $8$ hours). The newsvendor has limited transport capacities. Thus the total number of newspapers purchased in the morning is limited to a number of at most $v \in \N$ pieces. As usual for the newsvendor problem, there are costs $c_i \in \R_+$ per newspaper of type $i$ which has to be payed to the publisher. The selling price per newspaper of type $i$ is given a number $q_i \in \R_+$. A newspaper which was not sold during the day can be given back to the publisher for a return price $r_i \in \R_+$ per piece. It makes sense to assume $r_i \leq  c_i < q_i$. The newsvendor not only wants to maximize her daily gain but also wants to have a good balance between gain and life time. This leads to a multi-objective setting with two conflicting objectives: (i) maximize the gain (or equivalently minimize the loss), (ii) minimize the working time.

The average time to sell one newspaper of type $i \in \Set{1,\dots,n}$ at day $\omega \in \Omega$ is denoted by $t_i(\omega)$.
Usually, the lower the demand, the more time is required to sell a copy.

This problem can be formulated as a recourse problem.
	\leqnomode
	\begin{gather*}\tag{RP$_\text{NVP}$}\label{spr2a}
	  \Min_{\substack{x \in \R^n \\ y^1,\dots,y^N \in \R^m}} 
	  \sum_{j=1}^n 
	  \begin{pmatrix}
	   	c_j-r_j \\
	 	0
	   \end{pmatrix}
	 x_j + \frac{1}{N} \sum_{i=1}^N \sum_{j=1}^n 
		   \begin{pmatrix}
		   	r_j-q_j \\
		 t_j(\omega_i)
	  \end{pmatrix}
	 y^i_j \\ \notag
	 \quad\text{s.t.} \left\{ 
	 \arraycolsep=1.0pt
	 \begin{array}{rl}
	   	     \sum_{j=1}^n x_j\leq v &\\
	   	       0 \leq y^i \leq x &,\; i \in \Set{1,\dots,N} \\
	 \end{array}\right.
	\end{gather*}
	\reqnomode

Let us now consider a small specific example to illustrate the concept of preference for flexibility. We assume that there are $n=2$ types of newspapers: the ``Jena Post" (JP) and the``Bozen Times" (BT). The prices for the newspapers are shown in Table \ref{tab:prices}. Their demand observed over two working days is presented in Table \ref{tab:demand}.	
\begin{table}[ht]
    \begin{minipage}[t]{0.48\textwidth}
        \centering
        \begin{tabular}{lcc}
        \toprule
         & JP & BT \\
        \midrule
        Purchase Price & $\geneuro 3.00$ & $\geneuro 2.00$ \\
        Selling Price & $\geneuro 5.50$ & $\geneuro 5.00$ \\
        Return Price & $\geneuro 1.00$ & $\geneuro 2.00$ \\
        \bottomrule
        \end{tabular}
        \caption{Prices of newspapers}
        \label{tab:prices}
    \end{minipage}
    \hfill
    \begin{minipage}[t]{0.48\textwidth}
        \centering
        \begin{tabular}{lrr}
        \toprule
                 & JP & BT \\
        \midrule
        Monday & $200$ & $150$ \\
        Tuesday & $200$ & $100$ \\
 &   & \\
        \bottomrule
        \end{tabular}
        \caption{Demand for newspapers}
        \label{tab:demand}
    \end{minipage}
\end{table}
In the case of high demand, BT is clearly the better choice. The profit margin is higher, and there is no risk if the newspapers are not sold. However, we see that the demand for JP is higher, which is an advantage of JP over BT.

The time to sell one newspaper is assumed to be
$$ t_i(\omega) \coloneqq 200 / d_i(\omega).$$
Furthermore, let us assume that the total number of newspapers purchased in the morning is at most $v=200$, due to limited transport capacities. 

	How many copies of each type should the newsvendor purchase on Wednesday morning based on these data and her preferences with respect to a balance between gain and life time?

	If the efficient outcomes of the multi-objective linear program \eqref{spr2a} are presented to the newsvendor, she could be interested in the Pareto efficient expected outcome of $\geneuro 250$ with a working time of $100$ minutes. This expected outcome is obtained by purchasing $100$ pieces of JP and $0$ pieces of BT, so the first-stage decision could be $x=(100,0)^T$. However, the first-stage decision $x=(100,100)$ also contains a gain of $\geneuro 250$ with a working time of $100$ minutes as a possible outcome. If the newsvendor does not change her preferences, both options are suitable in the same way. But otherwise, the choice of $x=(100,100)$ provides more flexibility to the newsvendor, as shown in Figure \ref{fig:1}.

		\begin{figure}[ht]
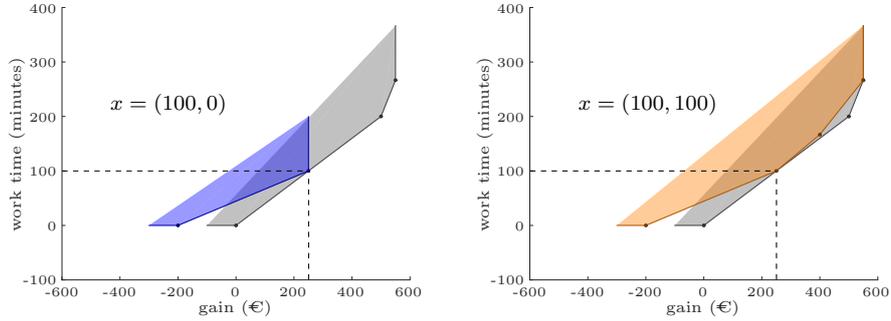

			\def\scalefactor{.4}
			\def\posx{100}
			\def\posy{200}
			\input{plot_single_trunc_modified_014.tex}
			\input{plot_rp2_trunc_modified_014.tex}
			\caption{A decision maker could be interested in the Pareto efficient expected outcome $\geneuro 250$ in $100$ minutes working time,
			i.e., in the point $y=(250,100)$. 
			There are several options to realize this outcome. One option is to purchase $100$ pieces of JP and no piece of BT, i.e., $x=(100,0)$.
			Another option is purchasing $100$ copies of JP and $100$ copies of BT, i.e., $x=(100,100)$.
			This ambiguity has consequences for the second-stage decision process whenever the newsvendor changes her preference.
			The grey sets in the background illustrate all possible outcomes (and outcomes worse than these) for all possible first-stage decisions $x$.
			The vertices and the edges between adjacent vertices represent the efficient outcomes.
			The colored sets display all possible second-stage outcomes (and outcomes worse than these) after fixing $x$. 
			We see that the efficient first-stage decision $x=(100,100)$ provides more flexibility in the second stage than the first-stage decision $x=(100,0)$.
			Moreover, $x=(100,100)$ is the most flexible among all first-stage decisions which enables the outcome of $\geneuro 250$ in $100$ minutes.
			Note that the $y_1$-axis is inverted as we maximize the gain while in the standard form of \eqref{spr1} we would minimize the loss.}\label{fig:1}
		\end{figure}

\end{ex}

\begin{rem}
\label{RemComparisonGH25}
Figure \ref{fig:1} also illustrates the difference between the solution concept in \cite{GuptaHunter25} and the one proposed here. According to \cite[Definition~1 and the subsequent paragraph]{GuptaHunter25}, the goal is to identify feasible first-stage decisions~$x$ whose expected outcomes are Pareto-minimal points within the set of all expected outcomes that can be obtained by some first-stage decision. Thus, both first stage decisions $x=(100,0)$ and $x=(100,100)$ would be optimal in the sense of \cite[Definition~1]{GuptaHunter25} showing that flexibility in the second stage is not taken into account.
\end{rem}

\section{Preferences for flexibility and set optimization} 
\label{SecSOandPfF}

Multi-objective stochastic recourse problems can be viewed as particular instances of decision making under a `preference for flexibility' \cite{Kreps79Economet}. The basic idea can be described as follows. The decision maker can select one among a number of sets of alternatives (called `menus' in \cite{Kreps79Economet}) and has to make the decision which alternative in the chosen set she
picks only at a later point in time (the `meal' from the chosen `menu'). This raises the question how to compare sets of alternatives instead of single alternatives.

If $\preceq$ is a preference relation on a set $\mathcal X$ assumed to be reflexive and transitive and minimization is a goal, then it is intuitive that a decision maker would prefer a set $A \subseteq \mathcal X$ over a set $B \subseteq \mathcal X$ if
\[
\forall b \in B, \exists a \in A \colon a \preceq b.
\]
In this case, we write $A \lel B$. A parallel relation with maximization in view is used in \cite{Kreps79Economet} where also completeness of $\preceq$ is assumed which is by no means necessary. The relation $\lel$ has the following property: if $A \supseteq B$, then $A \lel B$. This can be interpreted as a ` ``desire for flexibility"' \cite[p. 566 with formula (1.3)]{Kreps79Economet}: a decision maker prefers a situation with more options to choose from. The relation $\lel$ provides this feature.

Assume that a collection $\mathcal A \subseteq 2^{\mathcal X}$ of subsets of $\mathcal X$ is given and a decision maker has (or wants) to first choose one particular set out of $\mathcal A$ and then, in a later stage, choose an element out of the selected 
set. The information which might be interesting at the first stage is

$\bullet$ to know all her options which is the union $\bigcup_{A \in \mathcal A} A$

$\bullet$ to know minimal elements of $\mathcal A$ with respect to $\lel$, i.e., $\bar A \in \mathcal A$ satisfying
\[
A \in \mathcal A, A \lel \bar A \; \Rightarrow \; \bar A \lel A
\]
(there does not exist a set $A \in \mathcal A$ which is strictly better than $\bar A$ with respect to $\lel$; recall that minimization is the goal).

The relation $\lel$ on $2^{\mathcal X}$ is not antisymmetric in general even if $\preceq$ is on $\mathcal X$. However, one can restrict $\lel$ to a smaller set gaining antisymmetry but without loosing information. For $A \subseteq \mathcal X$, define
\[
\hat A = \{x \in \mathcal X \mid \exists a \in A \colon a \preceq x\}, 
\]
i.e., one adds all elements to $A$ which are worse then an element of $A$: 
as indicated above, this does no harm if minimization is the goal, since only alternatives are added that will not be chosen by the decision maker anyway. Then $\hat A \lel A \lel \hat A$ (a decision maker with preference $\lel$ would be indifferent between $A$ and $\hat A$) and $\hat A$ is the largest subset of $\mathcal X$ with this property. One can understand the operation to move from $A$ to $\hat A$ as a hull operator (the ``hat" operator). The set $\mathcal P(\mathcal X, \preceq) := \{A \in 2^{\mathcal X} \mid \hat A = A\}$ has very useful features which makes it accessible for many mathematical operations.

\begin{prop}
\label{PropFlexLattice}
One has $\hat A \supseteq \hat B$ if, and only if, $A \lel B$. The pair $\left(\mathcal P(\mathcal X, \preceq), \supseteq\right)$ is a complete lattice with 
\[
\inf_{A \in \mathcal A} A = \bigcup_{A \in \mathcal A} A
\]
for $\mathcal A \subseteq \mathcal P(\mathcal X, \preceq)$. In particular, $\lel$ is a partial order on $\mathcal P(\mathcal X, \preceq)$, i.e., also antisymmetric. 
\end{prop}

The proof is omitted; it can be found, e.g., for the case if $\preceq$ is a vector order in \cite[p. 68f]{HamelEtAl15Incoll}; the proof in the general case uses the same simple ideas.

The first part of the proposition indicates how the preference $\lel$ is transferred into a `preference for flexibility' \cite{Kreps79Economet}---
in this case directly into the $\supseteq$-relation. 

Moreover, Proposition \ref{PropFlexLattice} shows that the management goal of ``knowing all one's options'' is translated into the mathematical concept of an infimum in a particular complete lattice of sets which justifies the construction in \eqref{eq:upper_image_spr2}. On the other hand, the goal to find a ``best'' set becomes the problem of finding a minimal element with respect to $\supseteq$ in a collection of sets which are closed with respect to the hat operator: for such a minimal set, there does not exist a strictly larger feasible set, i.e., minimality with respect to $\supseteq$ means ``maximal flexibility" for the decision maker. Since $\supseteq$ is not a complete partial order on $\mathcal P(\mathcal X, \preceq)$ in general, there can be many minimal sets and the decision maker has to select one: this is the same dilemma one is faced with when many Pareto minimal alternatives exist.

The special case $\mathcal X = \mathbb \R^d$ with $\preceq$ being the component-wise order $\leq_{\R^d_+}$ 
generated by the convex cone $\R^d_+$ gives the link to the stochastic recourse problem, which is discussed in detail in the next section. 
In this specific setting, the hat operator is given by $\hat A = A + \R^d_+$ for  $A \subseteq \R^d$.
Let $F : \R^n \rightrightarrows \R^q$ be a polyhedral convex set-valued mapping with $F(x)= F(x)+\R^d_+$ for all $x \in \R^d$. Then the values $F(x)$ as well as the infimum over these values are convex polyhedra, see e.g., \cite[Proposition 3]{HeyLoe}, in particular, they are closed and convex sets. Thus we have $F(x) = \cl\conv F(x)$ for all $x \in \R^n$ and 
\begin{equation}
\label{EqFInfInG}
   \inf_{x \in \R^n} F(x) = \bigcup_{x \in \R^n} F(x) = \cl\conv \bigcup_{x \in \R^n} F(x).
\end{equation}
where $\cl S$ denotes the closure of a set $S \subseteq \R^d$, and the infimum is understood in $\{A \subseteq \R^d \mid A = \cl\conv(A + \R^d)\}$ with respect to $\supseteq$.
Therefore we define
\[
\mathcal G \coloneqq  \mathcal G(\R^d, \R^d_+) \coloneqq \Set{A \subseteq \R^d \mid A = \cl\conv(A + \R^d_+)},
\]
and the pair $(\mathcal G, \supseteq)$ will serve as image space for set-valued optimization problems occurring in this paper; it is a {\em complete lattice}, which means that every subset of $\mathcal G$ has an infimum and consequently also a supremum, see, e.g., \cite[Example 2.13]{HamelEtAl15Incoll} for more details. Infimum and supremum of a set $\mathcal A \subseteq \mathcal G$ with respect to $\supseteq$ are given by the following formulas:
\[
\inf_{A \in \mathcal A} A = \cl\conv\bigcup_{A \in \mathcal A} A, \quad \sup_{A \in \mathcal A} A = \bigcap_{A \in \mathcal A} A,
\]
where $\inf \emptyset = \emptyset$ and $\sup \emptyset = \R^d$ 
is understood.

Thus, it makes sense to ask for the infimum of a function $f:\R^n \to \mathcal G$ denoted
\[
\inf_{x \in \R^n} f(x).
\]
A set $\bar X \subseteq \R^n$ is called a \emph{solution} of this \emph{set minimization problem} if 

(S1) $\inf_{x \in \bar X} f(x) = \inf_{x \in \R^n} f(x)$,

(S2) each $\bar x \in \bar X$ satisfies
\[
x \in \R^n, \, f(x) \supseteq f(\bar x) \; \Rightarrow \; f(x) = f(\bar x).
\]
Condition (S1) can be considered as infimum attainment while (S2) states that each element of the solution set $\bar X$ is a minimizer of $f$ with respect to the order $\supseteq$ in $\mathcal G$. 
The upper image $\PP$ of the multi-objective linear program \eqref{spr2}, defined in \eqref{eq:upper_image_spr2}, can be seen as an infimum of the objective function
$f: S \to \mathcal G$, $f(z) = \Set{Pz} + \R^d_+$, that is, 
\[
 \PP = \inf_{z \in S} f(z).
\]
where the infimum is understood in $(\mathcal G, \supseteq)$. For more details on this solution concept for multiple objective linear programs the reader is referred to \cite{HeyLoe11, Loe11} and also \cite{HamelEtAl15Incoll}.

\section{Managerial decision making with preferences for flexibility}
\label{SecMDMwithPF}

Motivated by Example \ref{ex:1}, we now replace the multi-objective linear program \eqref{spr1} by a more general {\em polyhedral convex set optimization problem}. We will show that solving this set optimization problem furnishes the decision maker with additional information: It provides second-stage flexibility information for the finitely many first-stage decisions obtained from a solution to \eqref{spr1} and ensures maximal flexibility. The initial stochastic problem \eqref{spr} is replaced by
\leqnomode
\begin{gather*}\tag{RP$^\star$}\label{sprs}
 \text{minimize}\;\; E[\bm{Z}](x)\quad {s.t.} \quad 
 \left\{ 
   \begin{array}{l}
   	Ax=b,\; x \geq 0\\
 	x \in \R^n\text{,}
   \end{array}\right.
\end{gather*}
\reqnomode 
where $\bm{Z}:\R^n \to (2^{\R^d})^\Omega$ is the random polyhedral convex set-valued function 
$$\bm{Z}(x) \coloneqq \Set{C x + \bm{Q}\bm{y} \given \bm{T} x + \bm{W}\bm{y} = \bm{u},\; \bm{y}\geq 0  }.$$
The expectation $E[\bm{Z}]$ of $\bm{Z}$ is a deterministic polyhedral convex set-valued mapping.
The random polyhedral convex set-valued function $\bm{Z}$ can be expressed as
$$\bm{Z}(x) = \begin{pmatrix}
	Z(x,\omega_1)\\
	\vdots\\
	Z(x,\omega_N)
\end{pmatrix}$$
for the polyhedral convex set-valued mappings
$$ Z(\cdot,\omega_i):\R^n \rightrightarrows\R^d,\; Z(x,\omega_i)\coloneqq \Set{Cx+Q_i y \given T_i x + W_i y= u^i,\; y \geq 0}.$$ 
 The expectation of $\bm{Z}$ can be expressed using the notation of \eqref{spr2} as
\begin{gather*}
	 E[\bm{Z}]:\R^n \rightrightarrows \R^d,\\
	 E[\bm{Z}](x)\coloneqq \Set{P z \given \exists y^1,\dots,y^N \in \R^m:\; z=(x,y^1,\dots,y^N),\; z \in S}.
\end{gather*}
Problem \eqref{sprs} is a polyhedral convex set optimization problem. By
$$F \colon \R^n\rightrightarrows\R^d,\quad F(x)\coloneqq \left\{\begin{array}{cl}
	E[\bm{Z}](x) + \R^d_+ & \text{if } Ax = b,\; x\geq 0\\
	\emptyset    & \text{otherwise,}
\end{array}\right.$$
a polyhedral convex set-valued mapping $F$ is defined which actually maps into $\GG$.
Now, \eqref{sprs} can be expressed equivalently as
\leqnomode
\begin{gather*}\tag{RP$^{\star}$'}\label{sop}
 \text{minimize}\;\; F(x)\quad {s.t.} \quad 
   x \in \R^n \text{.}
\end{gather*}
\reqnomode
For a first-stage decision $x$ and a scenario $\omega \in \Omega$ we denote by $\PP(x,\omega) \subseteq \R^d$ the upper image of the second-stage decision problem \eqref{spr2xo}. 
\begin{prop}\label{prop:41n}
	For every first-stage decision $x$, the set $F(x) \subseteq \R^d$ is the expectation of the random polyhedron 
	$\PP(x,\cdot):\Omega \to 2^{\R^d}$.
\end{prop}
The proposition implies that the set optimization problem \eqref{sop}, which is just a reformulation of \eqref{sprs}, can be seen as the problem to minimize the expectation of the `optimal values' (upper images) of the second stage decision problems with respect to the first-stage decision variable $x$.

The upper image of \eqref{sprs} and \eqref{sop} is
$$ \PP^\star \coloneqq \bigcup_{x \in \R^n} F(x) + \R^d_+ = \bigcup_{x \in \R^n} F(x)\text{,}$$
where the last equation holds since the cone $\R^d_+$ was already added in the definition of $F$.
The upper image $\PP^\star$ can be seen as the optimal value of \eqref{sprs} and \eqref{sop} and is indeed the infimum of $F$ over $\R^n$ where the infimum is taken in $(\GG,\supseteq)$ as in \eqref{EqFInfInG}.

\begin{prop}\label{prop:33}
	The upper image $\PP$ of \eqref{spr}, \eqref{spr1} and \eqref{spr2} coincides with the upper image $\PP^\star$ of \eqref{sprs} and \eqref{sop}.
\end{prop}

We still assume that \eqref{eq:bounded} in Assumption (A2) holds, i.e., Problem \eqref{sop} is bounded. 
The solution concept introduced in the previous section via conditions (S1), (S2) can be specified to polyhedral problems as follows.

\begin{defn}\label{def:sol}
A \emph{solution} to the polyhedral convex set optimization problem \eqref{sop} is a finite set $\bar X \coloneqq \Set{\bar x^1,\dots,\bar x^k} \subseteq \R^n$ 
of points $\bar x \in \R^n$ satisfying 
\begin{equation}
\label{minimizer_sop}
	\not\exists x \in \R^n: F(x) \supsetneq F(\bar x) 
\end{equation}
with the additional property that each vertex of $\PP$ is contained in (at least) one the sets $F(\bar x^1),\dots,F(\bar x^k)$. 
\end{defn}

Note that $\supseteq$ has the meaning of `less than or equal to'.
Thus a point $\bar x\in \R^n$ satisfying \eqref{minimizer_sop} is called a {\em minimizer} of \eqref{sop}, while the latter property can be interpreted as {\em infimum attainment} and can be expressed equivalently as
\begin{equation}\label{inf_att}
	\PP = \conv \bigcup_{\bar x \in \bar X} F(\bar x) + \R^d_+
\end{equation} 
as this set is already closed and thus the closure in \eqref{EqFInfInG} can be omitted.
Infimum attainment corresponds to (S1) in Section \ref{SecSOandPfF}, but is expressed here for the specific setting used in this section. In addition to Section \ref{SecSOandPfF}, we assume here that the infimum is attained in a finite set $\bar X$. Condition \eqref{minimizer_sop} corresponds to (S2) in Section \ref{SecSOandPfF}.

Figure \ref{fig:1} illustrates this solution concept by the newsvendor problem from Example \ref{ex:1}.

	\begin{figure}[t]
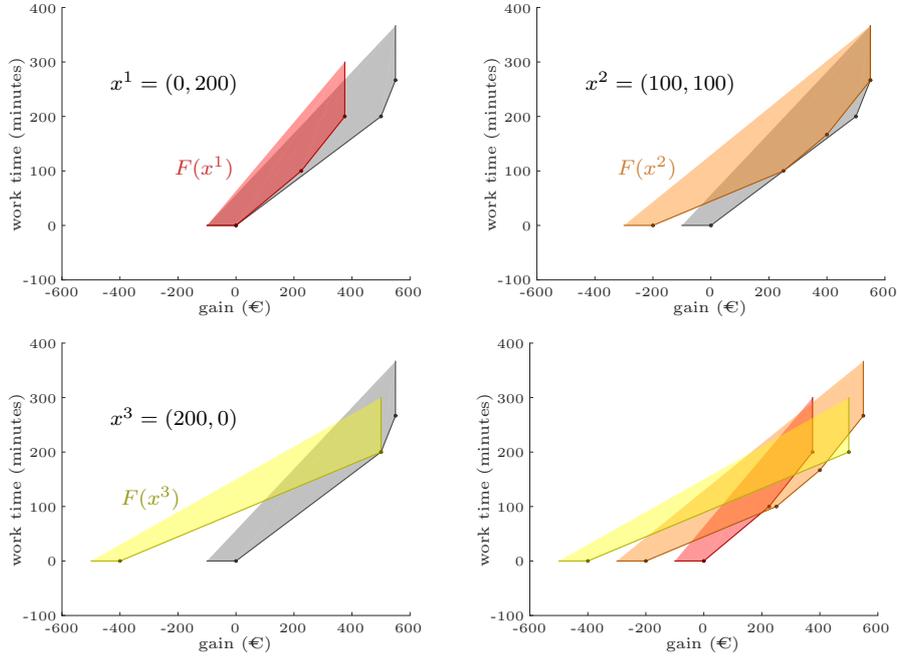

		\def\scalefactor{.4}
		\input{plot_rp1_trunc_modified_014.tex}
		\input{plot_rp2_trunc_modified_2_014.tex}
		\input{plot_rp3_trunc_modified_014.tex}
		\input{plot_rp_trunc_014.tex}
		\caption{Solution of the recourse problem from Example \ref{ex:1}. 
		The gray set in the background is the upper image $\PP = \PP^\star$.
		The colored sets in the front represent the objective values $F(x^1)$, $F(x^2)$ and $F(x^3)$.
		Every vertex of $\PP$ is contained in one of the sets $F(x^1)$, $F(x^2)$ and $F(x^3)$ (infimum attainment).
		The three displayed sets $F(x^i)$ are not properly contained in some other set $F(x)$, $x \in \R^n$ (minimality).
		Thus $\bar X = \Set{x^1, x^2 , x^3}$ is a solution to the polyhedral convex set optimization problem \eqref{sprs}.
		The last picture shows a comparison of the three solution elements.
		$F(x^i)$ represents the second-stage flexibility obtained for the first-stage decision $x^i$.}\label{fig:2}
	\end{figure}

The {\em vectorial relaxation} of the polyhedral convex set optimization problem \eqref{sop} is the multi-objective linear program 
\leqnomode
\begin{gather*}\tag{VR}\label{vr}
\text{minimize}\;\; v \quad{ s.t. } \quad 
\left\{ 
  \begin{array}{l}
	 v \in F(x) \\
	 x \in \R^n, v \in \R^d.
  \end{array}\right. 
\end{gather*}
\reqnomode
The constraints arising from $v \in F(x)$ contain the additional variables $y^1,\allowbreak\dots,\allowbreak y^N \in \R^m$
due to the definition of $F$ and $ E[\bm{Z}]$. Moreover, using the substitution $v =Pz$, we see that \eqref{spr2} is just a vectorial relaxation of \eqref{sop}. Since \eqref{spr2} is just a reformulation of \eqref{spr} and \eqref{sop} is a reformulation of \eqref{sprs}, we also have that \eqref{spr} is a vectorial relaxation of \eqref{sprs}.

The decision process for \eqref{sprs}, however, differs from the decision process for \eqref{spr}, as already described in Example~\ref{ex:1}. We describe two variants here, a ``forward'' and a ``backward'' approach.

In the “forward” approach, the decision maker begins the decision process by comparing the finitely many sets $F(x^1), \dots, F(x^k)$, where the set $\bar X = \Set{x^1, x^2, \dots, x^k}$ is a solution to \eqref{sop} according to Definition~\ref{def:sol}; see also Figure~\ref{fig:2} for an illustration. Recall that \eqref{sop} was derived from the original problem \eqref{spr}.
A suitable solution method can be found, for instance, in \cite{Loehne25SIOPT}. In the second stage, the multi-objective linear program \eqref{spr2xo} is solved for $x = \bar x$ and the realized scenario $\omega \in \Omega$. 

\subsubsection*{Decision Procedure 1 (DP1).}
\begin{description}
\item[Stage 1:] Choose $\bar x \in \bar X$, where $\bar X = \Set{x^1, x^2,\dots,x^k}$ is a solution to \eqref{sop}.
\item[Stage 2:] Choose a Pareto-minimal point of problem \eqref{spr2xo} for $x = \bar x$.
\end{description}

First, we remark that a solution $\bar X$ to the polyhedral convex set optimization problem is not necessarily unique, and there can even be infinitely many such solutions. An alternative solution provides a different form of flexibility. Moreover, the minimizers $x$ computed in this procedure enjoy the property that $F(x)$ contains a minimal point of $\PP$. However, problem \eqref{sop} may also have minimizers $x$ for which $F(x)$ does not contain a minimal point of $\PP$. An example, due to Frank Heyde, can be found in \cite[Example~3.6]{HamSch14}. Such minimizers may still be of interest to the decision maker, because they may offer a higher degree of flexibility in the second stage. This motivates the following ``backward" approach, which was already used in Example \ref{ex:1}, see Figure \ref{fig:1}.

\subsubsection*{Decision Procedure 2 (DP2).} 
\begin{description}
\item[Stage 1:] Choose one or several desired expected outcomes $z^1, \ldots, z^t \in \PP$ which are simultaneously realizable, i.e., there must be a first stage decision $x'$ such that $F(x')$ contains all the desired expected outcomes. In case of several desired expected outcomes, this typically leads to a trial and error process since such $x'$ does not need to exist in which case one has to modify the selected outcomes.
 Compute a minimizer $\bar x$ for \eqref{sop} such that $F(x') \subseteq F(\bar x)$.
\item[Stage 2:] Choose a Pareto-minimal point of problem \eqref{spr2xo} for $x = \bar x$.
\end{description}
The idea used in the first stage is extended in \cite{Loe24}, where a well-defined procedure for selecting desired expected outcomes is proposed, which yields a minimizer after finitely many steps.

\begin{rem} 
The difference between the two decision making procedures is that in (DP1) the first-stage variable $x$ is chosen without paying too much attention to the values $F(x)$ while (DP2) starts with a focus on potential outcomes of the second-stage decision. The latter might be especially useful if the availability of such outcomes is known beforehand. For example, a bank knows which portfolios are ``in the vault'' and can be used for risk compensation.
\end{rem}

\begin{rem} Clearly, it is assumed that the preference of the decision maker does not change from first to second stage decision. One could also consider problems where the preferences of the two stages are not aligned. This leads to completely different and even more complex problems.
\end{rem}

\section{The wait-and-see problem}
\label{SecWSP}

If both decisions, the choice of $x$ and $y$, are made after a scenario $\omega \in \Omega$ has been realized, one speaks about the {\em wait-and-see} approach. If the deterministic variable $x$ in \eqref{spr} is replaced by a random variable, we obtain the following {\em wait-and-see problem}.
\leqnomode
\begin{gather*}\tag{WS}\label{ws}
  \Min_{\bm{x} \in (\R^n)^\Omega, \bm{y} \in (\R^m)^\Omega} E[C \bm{x} + \bm{Q}\bm{y}]\quad\text{s.t.} \left\{ \begin{array}{rcl}
  	A \bm{x}&=&b\\
  	\bm{T} \bm{x} + \bm{W}\bm{y} &=& \bm{u}\\
  	\bm{x}, \bm{y} &\geq& 0\text{.}
  \end{array}\right. 
\end{gather*}
\reqnomode
Since \eqref{spr} arises from \eqref{ws} by the additional constraints $x(\omega_1) = \dots = x(\omega_N)$, the inclusion 
\begin{equation}\label{ineq1}
	\PP^\eqref{ws} \supseteq \PP^\eqref{spr}.
\end{equation}
holds for the upper images of \eqref{ws} and \eqref{spr}. In the single-objective case this inclusion reduces to the inequality $p^\eqref{ws} \leq p^\eqref{spr}$ \cite{Madansky60} for the optimal values of both scalar optimization problems. The difference $p^\eqref{spr} - p^\eqref{ws} \geq 0$ is called {\em expected value of perfect information}, see \cite{AvrielWilliams70OR} and the references therein. The additional constraints in the wait-and-see problem allow scenario-wise computation, which is useful if the \eqref{spr} instance is too large to be solved in practice. 

We now introduce a multi-objective counterpart. For each choice $v \in \PP^\eqref{spr}$ of the decision maker, the set 
$$\EVPI(v) \coloneqq \Set{x \in \R^d_+ \given v-x \in \PP^\eqref{ws}}$$
is called the {\em expected value of perfect information} with respect to the decision $v$. It consists of all possible ``improvements'' of $v$ if the scenario to be realized would be known. By \eqref{ineq1}, $\EVPI(v)$ is nonempty for all possible decisions $v \in \PP^\eqref{spr}$.

For $N\geq 2$, this inclusion can be strict as we see in the example below.
The multi-objective linear program \eqref{ws} can be decomposed into the $N$ smaller multi-objective linear programs 
\leqnomode
\begin{gather*}\tag{WS($\omega$)}\label{wsi}
  \Min_{x \in \R^n, y \in \R^m} C x + Q_\omega y \quad\text{s.t.} \left\{ \begin{array}{rcl}
  	A x &=&b\\
  	T_\omega x + W_\omega y &=& u_\omega\\
  	x,\, y &\geq& 0\text{.}
  \end{array} \right. 
\end{gather*}
\reqnomode
and the upper image $\PP^\eqref{ws}$ is obtained by the upper images of the smaller problems \eqref{wsi} as
$$ \PP^\eqref{ws} = \sum_{\omega \in \Omega} p_\omega \PP^\eqref{wsi}.$$

\begin{ex} \label{ex:2} (Multi-objective newsvendor problem II)\\
This is a continuation of Example \ref{ex:1}. The data set is now extended by the demand of Wednesday, see Table \ref{tab:3}.
\begin{table}
        \centering
        \begin{tabular}{lrr}
        \toprule
                 & JP & BT \\
        \midrule
        Monday & $200$ & $150$ \\
        Tuesday & $200$ & $100$ \\
        Wednesday & $50$ & $220$ \\
        \bottomrule
        \end{tabular}
        \caption{Demand for newspapers, extended data set}
        \label{tab:3}	
\end{table}

	\begin{figure}
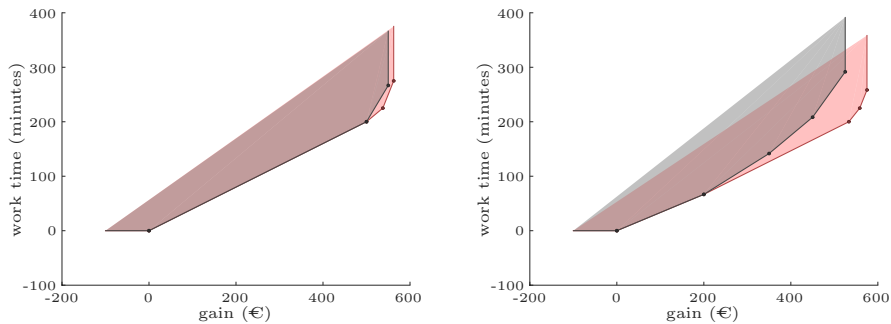

		\def\scalefactor{.4}
		\def\posx{100}
		\def\posy{200}
		\input{plot_cmp_ws_trunc_014.tex}
		\input{plot_cmp_ws_trunc_015.tex}
		\caption{Upper image of the wait-and-see problem (red, rear) compared to the upper image of the recourse problem (gray, front).
		Left: based on data from Monday and Tuesday only. Right: based on data from Monday to Wednesday.}\label{fig6}
	\end{figure}

The new data set is more volatile as the demand from Wednesday are quite different from the demand of Monday and Tuesday. In Figures \ref{fig6} we see the consequences for the wait-and-see problem. The more volatile data from Monday to Wednesday result in a larger ``difference'' between the upper images of the recourse problem and the wait-and-see problem, which can be interpreted as a ``higher'' value of information.
\end{ex}

\section{The expected-value problem}
\label{SecEVP}

The {\em expected value problem} arises from \eqref{spr} if the random variable $\bm{y}$ is replaced by a deterministic variable and all random data are replaced by their expectations. This leads to the multi-objective linear program 
\leqnomode
\begin{gather*}\tag{EV}\label{ev}
  \Min_{x \in \R^n, y \in \R^m} C x + E[\bm{Q}] y\quad\text{s.t.} \left\{ \begin{array}{rcl}
  	Ax&=&b\\
  	E[\bm{T}] x + E[\bm{W}] y  &=& E[\bm{u}]\\
  	x, y &\geq& 0\text{.}
  \end{array}\right. 
\end{gather*}
\reqnomode
Problem \eqref{ev} has fewer variables and fewer constraints than \eqref{spr}, in particular, this number does not depend on $N=|\Omega|$. 

\begin{prop}\label{prop:51}
	Under the assumption $\bm{Q}\equiv Q$ and $\bm{W}\equiv W$ the upper images of \eqref{ev} and \eqref{spr} satisfy the inclusion
	\begin{equation}\label{eq:52}
		\PP^\eqref{ev} \supseteq \PP^\eqref{spr}.
	\end{equation}
\end{prop}

In the single-objective case the inclusion \eqref{eq:52} reduces to the inequality $p^\eqref{ev} \leq p^\eqref{spr}$, see e.g.\ \cite{Birge82MP}, for the optimal values of corresponding scalar optimization problems.

A common approach in single-objective stochastic programming with recourse is to use the easier problem \eqref{ev} instead of \eqref{spr} in the first stage of the decision process. In the multi-objective setting this leads to the following polyhedral convex set optimization problem. 
\leqnomode
\begin{gather*}\tag{EV$^\star$}\label{evs}
 \Min_{x \in \R^n} \bar Z (x)\quad {s.t.} \quad A x = b, x \geq 0.
\end{gather*}
\reqnomode 
where $\bar Z : \R^n \to 2^{\R^d}$ is the polyhedral convex set-valued function 
$$\bar Z(x) \coloneqq \Set{C x + E[\bm{Q}] y \given E[\bm{T}] x + E[\bm{W}] y = E[\bm{u}],\; y \geq 0  }.$$
For \eqref{evs}, the first-stage decision process works completely analogous as for the recourse problem \eqref{sprs}. Eventually, the decision maker has chosen some minimizer $x \in \R^n$ of \eqref{evs}. 

After some scenario $\omega \in \Omega$ has been realized, the second-stage decision can be made by solving the multi-objective linear program \eqref{spr2xo}.
The decision maker chooses some minimizer $\bar y \in \R^m$ of \eqref{spr2xo}. The upper image $\PP^\eqref{spr2xo}$ is the set of all possible outcomes, together with all worse points, when scenario $\omega$ has been realized after the first-stage decision $x$. The expectation of the random set $\omega \mapsto \PP^\eqref{spr2xo}$ is the upper image of the multi-objective linear program
\leqnomode
\begin{gather*}\tag{EEV($x$)}\label{eevx}
  \Min_{\bm{y} \in (\R^m)^\Omega} E[C x + \bm{Q}\bm{y}]\quad\text{s.t.} \left\{ \begin{array}{rcl}
  	\bm{W}\bm{y} &=& \bm{u}- \bm{T} x\\
  	\bm{y} &\geq& 0\text{.}
  \end{array}\right. 
\end{gather*}
\reqnomode
which arises from the recourse problem \eqref{spr} by fixing $x$. The upper image of \eqref{eevx} can be expressed as
$$ \PP^\eqref{eevx} = \sum_{\omega \in \Omega} p_\omega \PP^{\eqref{spr2xo}} = E[\bm{Z}](x) + \R^d_+.$$
For feasible $x$, the latter expression is the value $F(x)$ of the set-valued objective function $F$ of the recourse problem \eqref{sprs}. 
This yields immediately the following inclusion.

\begin{prop}\label{prop:53} For arbitrary $x \geq 0$ with $A x = b$,
	\begin{equation}\label{eq:54}
		\PP^\eqref{spr} \supseteq \PP^\eqref{eevx}.
	\end{equation}
\end{prop}

Note that $\PP^\eqref{eevx}$ can be empty, i.e., (EEV($x$)) can be infeasible for a first-stage decision $x$ based on \eqref{evs}. In the single-objective case the inclusion \eqref{eq:54} reduces to the inequality $p^\eqref{spr} \leq p^\eqref{eevx}$, see e.g.\ \cite{Birge82MP}, for the optimal values of the associated scalar problems. 
We close this section by an example. 

\begin{ex} \label{ex:3} (Multi-objective newsvendor problem III)\\
The problem from Examples \ref{ex:1} and \ref{ex:2} can be used to show that the inclusion \eqref{eq:52} can be violated, see Figure \ref{fig8}. The reason is that $\bm{Q}$ is not constant.	
 	
	\begin{figure}[ht]
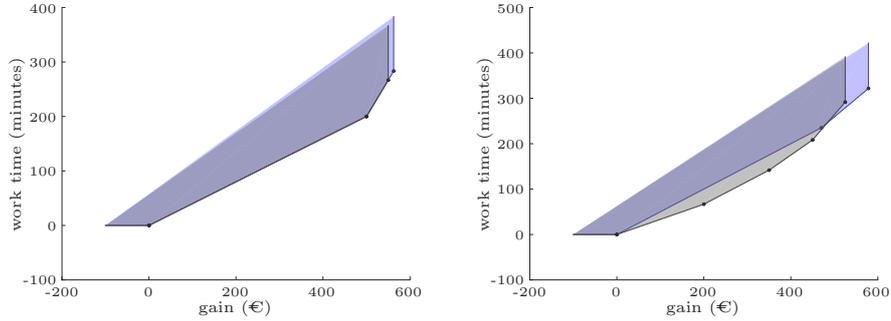

		\def\scalefactor{.4}
		\def\posx{100}
		\def\posy{200}
		\input{plot_cmp_ev_trunc_014.tex}
		\input{plot_cmp_ev_trunc_015.tex}
		\caption{Upper image of the expected value problem (blue) compared to the upper image of the recourse problem (gray). 
		Left: based on data from Monday and Tuesday only. Right: based on data from Monday to Wednesday, where we see that inclusion \eqref{eq:52} can is violated.}\label{fig8}
	\end{figure}
 
\end{ex}

\begin{ex} \label{ex:4} (Multi-objective newsvendor problem IV)\\
	Here we consider the problem from Example \ref{ex:1}. A solution of the expected value \eqref{ev} problem is
	$$ \bar X \coloneqq \left\{
	\begin{pmatrix}
		200\\
		0
	\end{pmatrix},\; 
	\begin{pmatrix}
		0 \\
		200
	\end{pmatrix},\; 
	\begin{pmatrix}
		75 \\
		125
	\end{pmatrix}\right\}.$$
	In Figure \ref{fig10} we see the upper images of \eqref{eevx}, where $x$ runs over $\bar X$. In particular, we see that inclusion \eqref{eq:54} is satisfied.

		\begin{figure}[ht]
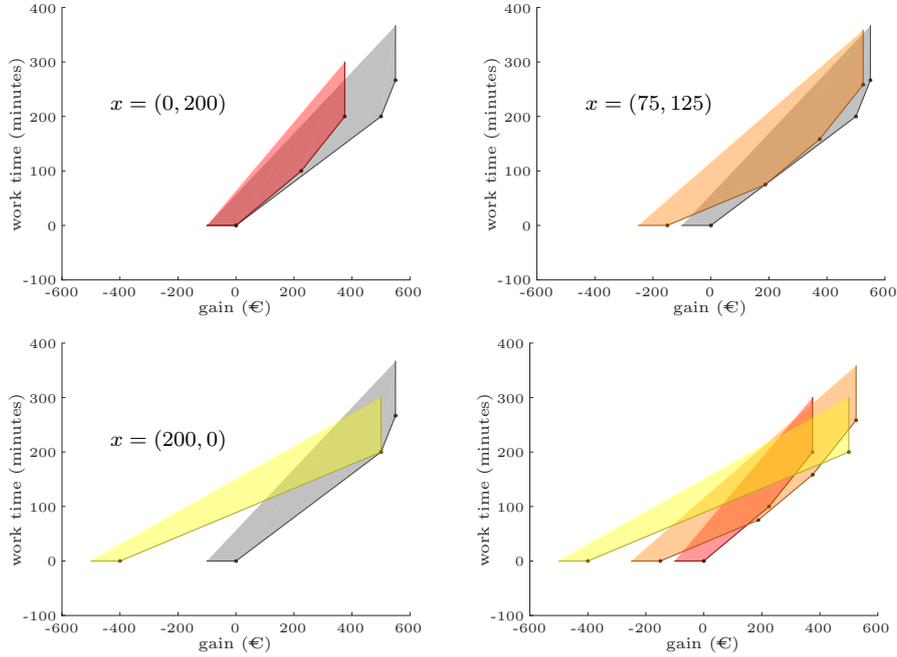

			\def\scalefactor{.4}
			\def\posx{100}
			\def\posy{200}
			\input{plot_eev1_trunc_014.tex}
			\input{plot_eev2_trunc_014.tex}
			\input{plot_eev3_trunc_014.tex}
			\input{plot_eev_trunc_014.tex}
			\caption{The upper images of \eqref{eevx} are displayed in the first three pictures. 
			The gray set in the background is the upper image of the recourse problem.
			We see that not all possible outcomes of the recourse problem can be obtained by this solution of \eqref{ev}, not even if we take the convex hull over the upper images of \eqref{eevx}.
			For instance, the maximal gain is not reached.}	\label{fig10}
		\end{figure}
	 
\end{ex}

\section{An application to risk management}
\label{section:ex_risk}

Let $\Omega$, $p$ be as in Section \ref{SecLRP}. A random vector $X \in (\R^{D \times 1})^\Omega$ with a natural number $D \geq 1$ is given which is now interpreted as the multidimensional random future payoff of a portfolio with $D$ assets. It is not assumed that all components of $X$ are denoted in the same reference instrument, e.g., USD. On the contrary, transaction costs may be present which makes it appropriate to consider set-valued risk measures. More background can be found in \cite{HamelHeyde21, HamHeyRud11}. The following assumptions are in force.




\begin{itemize}
\item There are $n$ eligible positions $H^1, \ldots, H^n \in (\R^{D \times 1})^\Omega$ in the market; they can be used for risk compensation and also serve as ``currencies" in the sense that a liquidation of a portfolio $X$ at the end of the trading/holding period is done into a linear combination of these positions. For example, all components of $H^i$ could be zero but one which is a low risk bond; the traditional interpretation is that the $H^i$'s non-zero components are currencies.

\item A solvency cone $\R^n_+ \subseteq K_0 \subseteq \R^n$ at time 0 serves as market model, i.e., proportional transaction costs are present for the eligible assets; $K_0$ can be generated by bid-ask spreads and is assumed to be a polyhedral cone facilitating exchanges of eligible positions into one another.
\item A random solvency cone $K_1 \in \left(2^{(\R^D)}\right)^\Omega$ with $\R^D_+ \subseteq K_1(\omega) \subseteq \R^D$ for all $\omega \in \Omega$ at terminal time is given; $K_1$ is assumed to be a random polyhedral cone again generated by (random) bid-ask spreads.
\end{itemize}
The following time consistency assumption is needed (and makes sense):
\[
x \in K_0 \quad \text{implies} \quad \sum_{i=1}^n x_iH^i \in K_1.
\]

\noindent
{\bf Trading.} A bank/an investor has an initial endowment $\bar x \in \R^n$ in units of the $n$ eligible portfolios; the simplest case is $n = 1$ with money in one currency; however, many deals involve the exchange of several assets such as stocks, physical assets like real estate and infrastructure, rights for exploitation of resources or patents etc. at the same time. The first transaction is to buy a portfolio $X$ available at the market at the price $\pi_0(X) \in \R^n$ using the initial endowment. Thus, the bank/investor can pick $X$ satisfying
\[
\pi_0(X) \in \bar x - K_0.
\]
For the sake of simplicity, the price operator $\pi_0 \colon (\R^{D \times 1})^\Omega \to \R^n$ is assumed to be linear. 

\medskip\noindent
{\bf Risk compensation.} The bank/investor wants/needs to compensate for the risk of $X$ by giving a deposit in $n$ eligible positions at time $0$. This deposit has the structure $\sum_{i=1}^n x_i H^i$ where $x_i$ is the number of units of position $i$ with uncertain future value $H^i$. The risk is evaluated by a risk measure $R$ given by
\[
R(X) = \left\{x \in \R^n \mid X + \sum_{i=1}^n x_iH^i \in K_1\right\}.
\]
The time consistency assumption ensures $R(X) + K_0 \subseteq R(X)$. The interpretation of the values $R(X)$ is as follows: $x \in R(X)$ is a risk management strategy where $x_i$ units of portfolio $H^i$ are given as deposit at initial time into which might include defaultable assets but are considered to be much more conservative, i.e., less risky investments than $X$. The overall portfolio $X + \sum_{i =1}^n x_iH^i$, i.e., $X$ plus a risk management strategy, should end up in a solvent position. This risk measure is a version of a worst-case risk measure and can of course be replaced by other types. 

In the finite setting used here, the graph of $R$ can be identified with a polyhedral set in $\R^{DN \times n}$.

\medskip\noindent
{\bf Objectives.} The position $\hat x = \bar x + x$ is available at initial time which is split into $\bar x$ used for investement and $x$ necessary for risk compensation. First, the bank/investor does not want to deposit more than necessary, i.e., $x \in R(X)$ should be minimal. Secondly, the gain from investing in $X$ should of course be maximal. At terminal time the position $X$ plus the deposit will be liquidated into the eligible positions according to the market at terminal time which is given by the random polyhedral convex cone $K_1$:
\[
\sum_{i=1}^n \bm{y_i} H^i \in X + \sum_{i=1}^n x_i H^i - K_1
\] 
with $n$ random coefficients $\bm{y_1}, \ldots, \bm{y_n} \geq 0$. The rational is that the investor gets back the deposit once she liquidates the position $X$ since no more deposit is necessary after liquidation. If a liquidation with only non-negative amounts $\bm{y_i}$ is not possible, the investor ends up in dept in at least one eligible position $H^i$ which is avoided by the non-negativity restriction.

Thus, the expression
\[
x - E[\bm{y}]
\]
needs to be minimized with respect to $\leq_{K_0}$ where the underlying assumption is that $\bm{y}$ is preferred over $\bm{y'}$ if $E[\bm{y - y'}] \in K_0$. 

\begin{rem}
One can understand the objective $x - E[\bm{y}]$ as putting the same weight on minimizing the deposit and maximizing the expected investment result. Different weights are possible such as $w_1x - w_2E[y]$ with $w_1, w_2 \geq 0$ and $w_1 + w_2 = 1$ (or even different weights for the components of $x$ and $E[\bm{y}]$) or one could formulate the problem as a vector optimization problem with objective $(x, -E[\bm{y}])$.
\end{rem}

The overall problem becomes
\[
\tag{RiskMin}
\text{minimize} \quad x - E[\bm{y}]
\]
\begin{align*}]
\text{subject to} \; & \; \pi_0(X) \in (\hat x - x) - K_0 \\
	 & \; X + \sum_{j =1}^n x_jH^j \in K_1 \\
	 & \; \sum_{i=1}^n \bm{y_i} H^i \in X + \sum_{i=1}^n x_i H^i - K_1 \\
	 & \bm{y_1}, \ldots, \bm{y_n} \geq 0 
\end{align*}
where $\bar x = \hat x - x$ is used in the first constraint. 

Problem (RiskMin) is a multi-objective stochastic linear program with recourse and $d=m=n$: the matrix $C$ in \eqref{spr} is the $n \times n$ unit matrix and $Q$ the negative of it (and thus is deterministic). Due to the assumption of polyhedrality for all involved cones the constraints can be transformed into linear equations which are deterministic (the first constraint) and stochastic (the second and the third) with suitable data $A$, $b$, $\bm{T}$, $\bm{W}$, $\bm{u}$. We omit the details at this point.

For this example, the ``flexibility function" $F$ is given by
\begin{multline*}
F(x) = E\left\{x - y(\omega) \,\big|\, X(\omega) + \sum_{i=1}^n x_iH^i(\omega)\in K_1(\omega), \right.\\
	\left.  X(\omega) + \sum_{i=1}^n x_iH^i(\omega) - \sum_{i=1}^n y_i(\omega)H^i(\omega) \in K_1(\omega), \, y_1, \ldots, y_n \geq 0  \right\} + K_0
\end{multline*}
if $ \hat x - \pi_0(X) - x \in K_0$ and $F(x) = \emptyset$ otherwise. In line with the desire for flexibility, the investor might want to have this set as large as possible.

An example illustrates the effects.

\begin{ex} 
Assume there are $D=3$ assets with $n=2$ of them being eligible.
Consider $N=2$ scenarios occurring with probabilities $p=(\tfrac{9}{10}, \tfrac{1}{10})$.
Furthermore, let the polyhedral convex cone $K_0$ be defined as the conic hull of the columns of the following matrices:
  \[
    \renewcommand{\arraystretch}{2}
    \begin{pmatrix} 1& -\frac{3}{10}\\ -\frac{3}{10}& 1\\\end{pmatrix}.\;
  \] 
One has $\R^2_+ \subseteq K_0$, so $K_0$ can be considered as a solvency cone for a two-asset-market. Similarly, define $K_1(\omega_1)$ and $K_1(\omega_2)$ via
  \[
    \renewcommand{\arraystretch}{2}
    \begin{pmatrix} 0& -\frac{7}{10}& 1& 1& 0& -\frac{4}{5}\\ -\frac{4}{5}& 0& -\frac{4}{5}& 0& 1& 1\\ 1& 1& 0& -\frac{7}{10}& -\frac{4}{5}& 0\end{pmatrix},\;
    \begin{pmatrix} 1& 0& -\frac{7}{10}& -\frac{4}{5}& 1& -0\\ -\frac{4}{5}& -\frac{7}{10}& 0& 1& 0& 1\\ 0& 1& 1& 0& -\frac{4}{5}& -\frac{4}{5}\end{pmatrix}.
\]
It also holds $\R^3_+ \subseteq K_1(\omega_1), K_1(\omega_2)$. We pick
\[
  H(\omega_1)=\renewcommand{\arraystretch}{2}\begin{pmatrix} 1& 1\\ 1& 1\\ 0& 0\\\end{pmatrix},\qquad
  H(\omega_2)=\renewcommand{\arraystretch}{2}\begin{pmatrix} 1& 0\\ 0& 1\\ 0& 0\\\end{pmatrix},
\]
which means that in scenario $\omega_1$ a portfolio with the first two components non-zero and equal is eligible while in scenario $\omega_2$ a portfolio is eligible if it only has a holding in asset 1 or in asset 2.

The time consistency condition can be checked using 
\[
x \in K_0 \quad \Leftrightarrow \quad x = \begin{pmatrix} s_1 - \frac{3}{10} s_2 \\ -\frac{3}{10}s_1 + s_2 \end{pmatrix}, \; s_1 \geq 0, s_2 \geq 0
\]
and
\[
x_1H^1(\omega_1) + x_2H^2(\omega_1) \in \R^2_+
\]
as well as
\[
x_1H^1(\omega_2) + x_2H^2(\omega_2) =
 \begin{pmatrix} s_1 - \frac{3}{10} s_2 \\ -\frac{3}{10}s_1 +  s_2 \\ 0 \end{pmatrix} \in K(\omega_2)
\]
if, and only if, 
\[
\begin{pmatrix} s_1 - \frac{3}{10} s_2 \\ -\frac{3}{10}s_1 + s_2 \\ 0 \end{pmatrix} =
\begin{pmatrix} 
	r_1 - \frac{7}{10}r_3 - \frac{4}{5} r_4 + r_5 \\ 
	-\frac{4}{5}r_1 - \frac{7}{10}r_2 + r_4 + r_6  \\ 
	r_2 + r_3 - \frac{4}{5}r_5 -  \frac{4}{5}r_6 
\end{pmatrix}, \;
r_1, \ldots, r_6 \geq 0.
\]
A solution of this system is given by $r_2 = \frac{4}{5}r_5$, $r_3 = \frac{4}{5}r_6$ and
\[
r_1 = \frac{983}{500}s_1, \; r_2 = \frac{69}{50}s_2, \; r_3 = \frac{69}{50}s_1, \; r_4 = \frac{983}{500}s_2,
\]
i.e., the time consistency condition is indeed satisfied.

Finally, let the random outcome be given by
  \[
    X(\omega_1)=\renewcommand{\arraystretch}{2}\begin{pmatrix} 3\\ 0\\ 1\\\end{pmatrix},\qquad
    X(\omega_2)=\renewcommand{\arraystretch}{2}\begin{pmatrix} -1\\ 1\\ -1\\\end{pmatrix}
  \]
and the initial wealth and price operator by
  \[
  \hat{x}=\renewcommand{\arraystretch}{2}\begin{pmatrix} 18\\ -1\\\end{pmatrix},\qquad 
  \pi_0=\renewcommand{\arraystretch}{2}\begin{pmatrix} 1\\ 1\\\end{pmatrix}.
  \]
Note that the second component of $\hat x$ is negative meaning that the agent is short in asset \#2 at initial time (but long in asset \#1). Buying back units of the second asset usually cannot be done at the same price as selling them which justifies the presence of the solvency cone $K_0$ which models the bid-ask spread of the market at initial time. See Figure~\ref{fig:11} for results with respect to this example. For the computations, the software package \texttt{bensolve tools} \cite{CirLoeWei17bt, CirLoeWei19} was used.

		\begin{figure}[ht]
			\def\scalefactor{.4}
			\def\posx{180}
			\def\posy{190}
			\input{plot_riskm_1_trunc-modified.tex}
			\input{plot_riskm_2_trunc-modified.tex}
			\caption{Both pictures show the set of all possible outcomes, i.e., the upper image of the problem, displayed in gray in the background. The outcomes $F(x^1)$ (left) and $F(x^2)$ (right), which can be obtained in the second stage after the first-stage decisions $x^1=(4.475,-1)$ and $x^2=(4.475,1.758)$ have been made, respectively, are shown in the foreground.
A multi-objective decision maker who prefers the payoff $p=(4.253,-5.728)$ (the highlighted white point) would be satisfied with both decisions $x^1$ and $x^2$, as both generate the point $p$ on the Pareto frontier.
However, a set-valued decision maker would not choose $x^1$, as it does not provide maximal flexibility. Indeed, $F(x^2)$ is a proper superset of $F(x^1)$. This means that the first-stage decision $x^2$ also allows obtaining the payoff $p$ in the second stage, while additionally offering greater flexibility. Here $x^2$ provides even maximal flexibility.}\label{fig:11}
		\end{figure}
	 
\end{ex}

\section{Conclusions.}
\label{SecConclusions}

This paper is the first which uses tools from set optimization to construct management strategies for multi-objective stochastic problems with recourse. This new approach can be understood as a managerial implementation of preferences for flexibility known in economics. The flexibility goal could not be reached without the set-valued approach; it would remain invisible if the focus would only be on Pareto minimal points and/or (weighted sum) scalarizations.

Examples illustrate the new effect. In particular, a two-dimensional newsvendor problem is discussed to show the effects of the new solution approach. An application to risk management in the presence of transaction costs confirms applicability to finance.

The problem considered in this paper can be extended in several directions. 
Assumptions can be relaxed such as the boundedness assumption in (A2) and linearity. Note that the set-valued approach works well also for more general cones than just $\R^d_+$ as in (A1).
The model can be extended inductively to problems with more than two stages. It would also be interesting to study more complex applications with sustainability objectives in addition to traditional management goals such as cost minimization or profit maximization. For risk management models, more complex risk measures as well as more sophisticated market models can be considered. A version for nonatomic probability spaces, with a different solution concept, is discussed in \cite{GuptaHunter25} which leads to a problem that is, in general, infinite-dimensional.

Bensolve tools \cite{CirLoeWei17bt} were used for the computations in this paper. However, the development of numerical methods for (polyhedral convex) set optimization problems is still at its beginning \cite{Loehne25SIOPT}. To reduce the computational effort, a scenario reduction method \cite{DupavcovaEtAl03} applied to problem \eqref{spr1} might be useful (see also \cite{BertsimasEtAl23MS}).

\section{Appendix}

\textbf{Proof of Proposition \ref{prop:31}.} Consider an arbitrary $i \in \Set{1,\dots,N}$. Using the notation $Q_i = Q(\omega_i)$, $T_i = T(\omega_i)$, etc., we observe that $\bar y^i$ is feasible for \rm{(RP$_2(\bar x,\omega_i)$)}. Suppose $\bar y^i$ is not optimal. Then there is some $\hat y^i$ that is feasible for \rm{(RP$_2(\bar x,\omega_i)$)} with $Q_i \hat y^i < Q_i \bar y^i$. 
If we replace  $\bar y^i$ by $\hat y^i$ in the solution $(\bar x,\bar y^1,\dots, \bar y^N)$, this point remains feasible for \eqref{spr1}. Because $p_i$ is assumed to be positive, this replacement results in a smaller objective value. This contradicts the assumption that $(\bar x,\bar y^1,\dots, \bar y^N)$ is a minimal solution. 
Thus $\bar y^i$ must be optimal for all $i \in \Set{1,\dots,N}$.\medskip

\noindent
\textbf{Proof of Proposition \ref{prop:41n}.} The upper image of \rm{(RP$_2(\bar x,\omega_i)$)} can be expressed as
$$ \PP(x,\omega_i) = Z(x,\omega_i) + \R^d_+.$$
So, the expectation of the random set $\PP(x,\cdot):\Omega \to 2^{\R^d}$ is
$$
	 E[\PP(x,\cdot)] = E[\mathbf{Z}+\R^d_+] = E[\mathbf{Z}] + \R^d_+.
$$	  
Because a first-stage decision $x$ is feasible, i.e., $Ax=b$ and $x \geq 0$, we have $E[\PP(x,\cdot)]=F(x)$.\medskip

\noindent
\textbf{Proof of Proposition \ref{prop:33}.} Since \eqref{spr}, \eqref{spr1}, and \eqref{spr2} are merely reformulations of the same problem, they all have the same upper images. The same holds for
\eqref{sprs} and \eqref{sop}. Thus it remains to show that \eqref{spr2} and  \eqref{sop} have the same upper images. This means we need to show that
$$ \Set{Pz \given z \in S}+\R^d_+ = \bigcup_{x \in \R^n} F(x).$$
For feasible $x$, i.e., $Ax=b$ and $x\geq 0$, we have
\begin{align*} F(x) &= E[\mathbf{Z}](x)+\R^d_+ \\
&=\Set{Pz \given \exists y^1,\dots,y^N \in \R^m:\; z=(x,y^1,\dots,y^N),\; z \in S} + \R^d_+,
\end{align*}
and otherwise, $F(x)=\emptyset$.
Thus, the union over all $F(x)$ results in $\Set{Pz \given z \in S}+\R^d_+$.

\medskip

\noindent
\textbf{Proof of Proposition \ref{prop:51}.} The result is derived from the following two facts: (i) Feasibility of $(x, y^1,\dots y^N)$ for \eqref{spr} implies feasibility of $(x,\sum_{i=1}^N p_i y^i)$ for \eqref{ev} if $\bm{W}\equiv W$; (ii) The objective function of \eqref{spr} at $(x, y^1,\dots y^N)$ has the same value as the objective function of \eqref{ev} at $(x,\sum_{i=1}^N p_i y^i)$ if $\bm{Q}\equiv Q$.

\bibliography{ref-stoch-recourse}

\end{document}